\documentclass{amsart}
\usepackage[all]{xy}
\SelectTips{cm}{}
\usepackage{amsmath}
\usepackage{amssymb}
\usepackage{amscd}
\usepackage{amsthm}
\usepackage{amsfonts}
\usepackage{amsxtra}
\theoremstyle{plain}

\newtheorem{bigthm}{Theorem}
\theoremstyle{definition}

\newtheorem*{Rem-intro}{Remark}

\newcommand{\Dirlim}{\varinjlim}

\newcommand{\ZZ}{{\mathbb{Z}}}
\newcommand{\QQ}{{\mathbb{Q}}}

\newcommand{\CC}{{\mathbb{C}}}

\newcommand{\idc}{_{\circ}}

\newcommand{\KK}{{\mathcal{K}}}
\def\:{\colon\!}

\begin{document}

\title[From Rational Homotopy to K-Theory]
{From Rational Homotopy to K-Theory for Continuous Trace Algebras}
 \author[Klein]{John R. Klein}
 \address{Department of Mathematics,
      Wayne State University,
   Detroit MI 48202}
 \email{klein@math.wayne.edu}
\author[Schochet]{Claude L.~Schochet}
\address{Department of Mathematics,
     Wayne State University,
     Detroit MI 48202}
\email{claude@math.wayne.edu}
\author[Smith]{Samuel~B.~Smith}
\address{Department of Mathematics,a
     Saint Joseph's University,
     Philadelphia PA 19131}
\email{smith@sju.edu}

\begin{abstract}
Let $A$ be a unital $C^*$-algebra. Its unitary group, $UA$, contains a wealth of topological
information about $A$. However, the homotopy type of $UA$ is out of reach even for $A = M_2(\CC)$.
There are two simplifications which have been considered. The first, well-traveled road, is to pass
to $\pi_*(U(A\otimes \KK ))$ which is isomorphic (with a degree shift) to $K_*(A)$. This approach has led
to spectacular success in many arenas, as is well-known.

A different approach is to consider $\pi _*(UA)\otimes\QQ $, the rational homotopy of $UA$.  In joint work with G. Lupton  and N. C. Phillips we have calculated
 this functor for the cases $A = C(X)\otimes M_n(\CC )$ and $A$ a unital
continuous trace $C^*$-algebra.   In this note we look at some concrete examples of this calculation and, in
particular, at the $\ZZ$-graded map
\[
\pi _*(UA)\otimes\QQ \longrightarrow K_{*+1}(A)\otimes\QQ  .
\]

\end{abstract}
\maketitle
%\centerline{  11-4-08}
\setcounter{tocdepth}{1}
%\tableofcontents

\section{Statement of the Main Theorem }

  Let $M_n = M_n(\CC)$ be the complex matrices, $U_n$  its group of unitaries, and let $PU_n$ be
  the quotient of $U_n$ by its center.
Let $\zeta \: T \to X$ be a principal $PU_n$-bundle over a compact space $X$, let $PU_n$ act on
 $M_n  $ by conjugation and let
$$
 T \times_{PU_n} M_n  \to X
$$
be the associated $n$-dimensional complex matrix bundle. Define
 $A_\zeta$ to be the set of continuous sections of the latter.
These sections
have natural pointwise addition, multiplication, and $*$-operations
and give $A_\zeta $ the structure of a  unital $C^*$-algebra. The algebra
$A_\zeta$ is  the most general
 unital  continuous trace $C^*$-algebra.
Let $UA_\zeta $ denote the topological group of unitaries of $A_\zeta $.
We have determined the rational homotopy type of $UA_\zeta.$

To state our calculation of the rational homotopy groups of $UA_\zeta$,
we introduce  some notation.   Given   graded $\mathbb Z$-graded vector spaces $V$ and
$W$, we grade the tensor product $V\otimes W$  by declaring that $v\otimes w$ has
grading $|v| + |w|$. Let $V\widetilde{\otimes} W$ be the effect of considering
only tensors with non-negative grading.

Given   elements $x_1, \ldots, x_n$, each of homogeneous degree,
write $\langle x_1, \ldots,  x_n \rangle$ for the graded vector space  with basis $x_1, \ldots, x_n$.
Given a topological group $G$,
write $G\idc$ for the path component of the identity in $G$.

The following is the principal result of our recent paper \cite{KSS}.

\begin{bigthm} \label{bigthm:mainA}\cite{KSS}
Let $\zeta$   be  a principal $PU_n$
 bundle over  a compact metric space $X$.  Let
$A_\zeta$ be the associated continuous trace $C^*$-algebra, and let $UA_\zeta$ its group of unitaries.
Then the rationalization of $(UA_\zeta)\idc$ is rationally $H$-equivalent to a product of
rational Eilenberg-Mac\,Lane spaces  with the standard multiplication,
with degrees and dimensions corresponding to an isomorphism of graded vector spaces
$$
\pi_*\left((UA_\zeta)\idc \right) \otimes \QQ \,\, \cong \,\, \Check{H}^*(X; \QQ) \,
\widetilde{\otimes}  \,  \langle s_1, s_3, \ldots, s_{2n-1}\rangle \, .
$$
\end{bigthm}

{\flushleft (}In the above,  $\Check{H}^*(X; \QQ)$ denotes the \v{C}ech cohomology of $X$
with rational coefficients and
we follow the convention that cohomology is graded in degrees $\le 0$. The basis element
$s_{2j-1}$ has degree $2j-1$.)
\medskip

   \section{Stabilization and the first example: $A = M_n$}

   We wish to explore the relationship between this result and the classical results on the $K$-theory
    of  continuous trace $C^\ast $-algebras. Recall first that homotopy and $K$-theory are related by stabilization.
  $K$-theory for  unital $C^*$-algebras is defined by taking $K_0(A)$ to be the Grothendieck group
of finitely generated projective left $A$-modules. For $j > 0$,
    write $U_\infty A =  \Dirlim U_nA $ with the weak topology and then define $K$-theory   for $j > 0$  by
    \[
    K_j(A) = \pi _{j-1}(U_\infty A) .\
    \]
  Thus there is a natural stabilization map
     \[
     \pi _j(UA) \longrightarrow \pi _j(U_\infty A) \cong K_{j+1}(A)
\]
which we denote
\[
\sigma ^{integral}  :  \pi _*(UA) \longrightarrow   K_{*+1}(A) .
\]
At this point the functor $K_j(A)$ is defined for $j \geq 0$.
  Bott periodicity implies that
\[
K_*(A) \cong K_{*+2}(A)
\]
for $* > 0$. So it is natural to regard $K_*(A)$ as a $\ZZ /2$-graded theory (and we confess to having
done so in the past.)
  However, for our present purpose it will be vital NOT to do so. That is,  {\emph{in this note we regard $K$-theory as defined for all non-negative integers.}}

 Ideally we would like to compute  the image of $\sigma ^{integral}$. However, this is probably as difficult
 as computing $   \pi _j(UA) $, which is out of reach even in the simplest cases. Hence we have been focusing on the rationalization
 of this group which, while far simpler and hence carrying less information, has the advantage of being computable.
 Thus we focus on the rational stabilization  map
 \[
\sigma :  \pi _*(UA)\otimes\QQ  \longrightarrow   K_{*+1}(A)\otimes\QQ .
\]
The first case to consider is $A = M_n$. Then we are looking at
 \[
\sigma :  \pi _*(U_n)\otimes\QQ  \longrightarrow   K_{*+1}(M_n )\otimes\QQ .
\]
Now
\[
\pi _*(U_n) \otimes\QQ \,\,\cong\,\, \langle s_1, s_3 , \dots, s_{2n - 1} \rangle
\]
with $|s_{2j-1} | = 2j-1$
and of course $\sigma $ is linear. So
$\sigma $ is determined by the various  $\sigma (s_{2j-1} ) $. We know by various easy arguments that
these map non-trivially, so that there is an isomorphism
\[
\sigma : \pi _{2j - 1}(U_n)\otimes\QQ  \overset\cong\longrightarrow  K_{2j}(M_n)\otimes\QQ
\]
for $j = 1, \dots n$ and hence the range of $\sigma $ is exactly
\[
\bigg\{ K_{2}(M_n)\otimes\QQ ,\,\,\,\,\, K_{4}(M_n)\otimes\QQ ,\,\,\,\dots ,\,\,\, K_{2n}(M_n)\otimes\QQ  \bigg\} .
\]

 In particular,   if we regard the range of $\sigma $ as an invariant of
$A$ with values in $\ZZ ^+$-graded $K$-theory, then this invariant distinguishes between $M_n $ and $M_k$ for $n \neq k$.

Note that we are NOT saying that the classes $\sigma (s_{2j-1})$ are
 actually in $K_*(M_n )$.
  Consider the diagram
\[
\begin{CD}
\pi _*(U_n) @>>>    \pi_*(U_n)\otimes\QQ @>\sigma >> K_{*+1}(M_n)\otimes\QQ  \\
@Vh'VV    @VhVV    \\
H_*(U_n ;\ZZ ) @>>>  H_*(U_n ; \QQ )    \\
\end{CD}
\]
where
 the vertical maps are the Hurewicz maps.
Filling in three well-understood groups, this diagram
becomes
\[
\begin{CD}
\pi _*(U_n) @>>>  \langle s_1, s_3, \dots , s_{2n-1}\rangle @>\sigma >> K_{*+1}(M_n)\otimes\QQ  \\
@Vh'VV    @VhVV    \\
\Lambda _\ZZ (g_1, g_3, \dots g_{2n-1} )  @>>>  \Lambda _\QQ (g_1, g_3, \dots g_{2n-1} )    \\
\end{CD}
\]
where
 $\Lambda _R$ denotes the exterior algebra over the ring $R$ with given generators
and
  $h(s_j) = g_j$ for each $j$.
The map $h' $ is not onto the generators.

\section{Second Example: $X = S^3$}

In order to make Theorem A concrete,   focus upon the special case where $X = S^3$.   Jonathan Rosenberg \cite{RosHom}  uses this case to illustrate beautifully the role of the Dixmier-Douady invariant in integral twisted
$K$-theory, and that paper should be consulted for the exact relationship, particularly how this plays out
in his twisted Atiyah-Hirzebruch spectral sequence.  

 Theorem A
asserts
  that the space of unitaries $UA_\zeta $ is rationally equivalent to the space of functions $F(S^3, U_n)$.
So let us examine   this space carefully.

Note first that its path components are interesting. Fix a base point $x_0 $ for $S^3$. There is a standard fibration
\[
F_\bullet (S^3, U_n) \longrightarrow F(S^3, U_n) \overset{p}\longrightarrow  U_n
\]
where $p(f) = f(x_o) $ and $F_\bullet $ denotes base point preserving maps. This fibration has a section (send
a point $u \in U_n $ to the constant map $S^3 \to U_n$ that takes every element of $S^3$ to $u$) and hence
there are short exact sequences in each degree
\[
0 \to \pi _*(F_\bullet (S^3, U_n))  \longrightarrow \pi _*(F(S^3, U_n))  \overset{p_*}\longrightarrow  \pi _*(U_n) \to 0.
\]
In particular, since $U_n$ is connected,
\[
\pi _0(\,F(S^3, U_n)\,) \cong \pi _0(F_\bullet (S^3, U_n)) \cong \pi _3(U_n) \cong \ZZ .
\]
The generator of $\pi _3(U_n) \cong \ZZ$ is given by the natural composition
\[
S^3 \cong SU_2 \to U_2 \to U_n
\]
where $U_2 $ is included in $U_n $ via $u \to u\oplus 1$.

In higher degrees after rationalization we obtain for each $k$ the short exact sequence
\[
0 \to \pi _k(F_\bullet (S^3, U_n))\otimes\QQ  \longrightarrow \pi _k(F(S^3, U_n))\otimes\QQ  \overset{p_*}\longrightarrow  \pi _k(U_n)\otimes\QQ \to 0
\]
This helps us understand the result for Theorem A, which states (in this case) that
\[
\pi _*(\, F(S^3, U_n)\,)\otimes\QQ  \cong H^*(S^3;\QQ )\otimes \langle s_1, s_3, \dots , s_{2n-1}\rangle .
\]
Write
\[
H^*(S^3; \QQ ) = \langle\,1,\,\, x_3 \,\rangle
\]
 with $x_3 $ denoting the generator in dimension $3$. Then
$\pi _*(F(S^3, U_n))\otimes\QQ $ is spanned by two types of classes. There are the classes $1\otimes s_{2j-1}$ of degree $2j-1$ and
the classes $x_3 \otimes s_{2j-1}$ of degree $2j-1 - 3 = 2j - 4$. The short exact sequence
\[
0 \to \pi _*(F_\bullet (S^3, U_n))\otimes\QQ  \longrightarrow \pi _*(F(S^3, U_n))\otimes\QQ  \overset{p_*}\longrightarrow  \pi _*(U_n)\otimes\QQ \to 0
\]
(which splits, of course, since these are all rational vector spaces) becomes
\[
0 \to \langle x_3 \otimes s_3, \dots x_3\otimes s_{2n-1} \rangle \,\,\longrightarrow\,\, \pi _*(F(S^3, U_n))\otimes\QQ
\,\,\longrightarrow\,\, \langle 1\otimes s_1, \dots 1\otimes s_{2n-1}\rangle \to 0
\]
 Note that the class $x_3\otimes s_1$ is {\emph{not}} present since it would have negative degree.

Now, what happens when we map to $K$-theory?
First, the fact that $UA_\zeta $ is rationally equivalent to $F(S^3, U_n)$ implies that
\[
K_*(A_\zeta )\otimes\QQ \cong K_*(C(S^3))\otimes\QQ  \cong K^*(S^3)\otimes\QQ  \cong \QQ
\]
in every degree. The generator in even degree is simply the class of the trivial line bundle (i.e. to the one
dimensional trivial projection)  and the
generator in odd degree corresponds to the Bott generator in that degree.
An easy naturality argument using the result of the previous
section implies that the class $1\otimes s_{2j - 1} $ maps to the class in $K_{2j}(A_\zeta )\otimes\QQ $ that corresponds
to a multiple of the one-dimensional trivial projection in $K_0(A_\zeta )$ under the Bott map as is the case
for $A = M_n$.

The other classes are more interesting. The class $x_3\otimes s_{2j-1} $
has degree $2j - 4$ and hence must map to $K_{2j-3}(A_\zeta )$. The first example is $x_3 \otimes s_3$
mapping to $K_1(A_\zeta)\otimes\QQ $ and the last is $x_3\otimes s_{2n-1}$ mapping to $K_{2n-3}(A_\zeta)\otimes\QQ $.

To summarize: if $X = S^3$ then, independent of the bundle $\zeta $, the image of the stabilization map
\[
\sigma :   \pi _*(UA_\zeta ) \otimes\QQ \longrightarrow K_{*+1}(A_\zeta ) \otimes\QQ
\]
has basis elements
\[
\sigma (1 \otimes s_{2j - 1}) \in K_{2j}(A_\zeta )\otimes\QQ \qquad\qquad j = 1,\dots , n
\]
and
\[
\sigma (x_3 \otimes s_{2j-1}) \in K_{2j - 3}(A_\zeta )\otimes\QQ \qquad\qquad j = 3, \dots n .
\]
  Thus the image of the
stabilization map $\sigma $ consists of the groups
\[
\bigg\{
K_0(A_\zeta )\otimes\QQ ,\,\,
K_1(A_\zeta )\otimes\QQ ,\,\,
\dots ,\,\,
 K_{2n-3}(A_\zeta )\otimes\QQ ,\,\,
K_{2n-2}(A_\zeta )\otimes\QQ ,\,\,
K_{2n}(A_\zeta )\otimes\QQ
\bigg\}
\]
and no others.

It is interesting to compare this result with the integral $K$-theory of $A_\zeta $. Rosenberg \cite{RosHom} 
shows that if the Dixmier-Douady invariant is non-zero then $K_*(A_\zeta )$ vanishes for $*$ even and is a 
finite cyclic group for $*$ odd. This implies that $K_*(A_\zeta ;\QQ) = 0$.   Thus there is no hope of trying 
to write the rationalization of the space $UA_\zeta $ when the Dixmier-Douady invariant is non-trivial as $U(A_\zeta \otimes N)$ for some $C^*$-algebra in the bootstrap
category \cite{Sch}  with $K_0(N) = \QQ $ and $K_1(N) = 0$.

\section{If one grading is good, then two gradings...}

We note that Theorem A gives a natural bigrading to $\pi _*(UA_\zeta )\otimes\QQ $.
In the $X = S^3$ example considered above,
 the classes $1\otimes s_{2j-1}$
have bidegree $(0,2j-1)$ and total degree $2j-1$, and the classes $x_3 \otimes s_{2j-1}$ have bidegree $(-3, 2j-1)$ and total degree
$-3+(2j-1) = 2j - 4$.
The bigrading has quite a bit of naturality associated with it. In the simplest case, with $A = C(X)\otimes M_n $
it is keeping track both of the cohomology degree and of the size of the matrix!

This bidegree is of course completely lost when passing to $\ZZ /2 $-graded $K$-theory. For an elementary
example, consider the case
 $X = \CC P^2$.
The rational (indeed, integral in this case) cohomology ring is a truncated polynomial
 algebra on a generator  $c \in H^2(CP^2)$
with $c^3 = 0$. Take $A =   F(\CC P^2, M_3)$. Then  the classes $c\otimes s_3 $ and $c^2\otimes s_5$
have different bidegrees, and hence are distinguished, but they have the same total degree and hence
have the same degree when one passes to $K$-theory, even when $K$-theory is $\ZZ ^+$-graded!

  Now, take $A = C(X)\otimes M_n$.  One might reasonably ask for a calculation of $[A, A]$,
  the homotopy classes of unital $*$-homomorphisms $[A, A]$.
  One source of such maps are the induced maps from functions $f: X \to X$ and so one might hope to determine
 $[A, A]$ as some functor of $[X, X]$. (Determining $[X,X]$ itself is extremely difficult even in fairly simple cases.).
 Another source of maps are the induced maps from $*$-homomorphisms $M_n \to M_n$. These are known, of course. Non-trivial maps must be isomorphisms, every isomorphism is inner, and hence every such map is given by conjugation by a unitary, so we are back to $PU_n$.  The real problem is that there are other maps besides these
 two types that intertwine the two.

 There is a   natural commuting  diagram
  \[
  \begin{CD}
  [A, A] @>\phi _K>>            End _{\ZZ /2}  (K^*(X))  \\
    @VV\phi V         @VVV    \\
    End _{**} \big( \pi _*(UA)\otimes \QQ \big)    @>>> End _{\ZZ /2}  ( K^*(X)\otimes\QQ
  \end{CD}
  \]
It might seem at first glance that one would be better off using $\phi _K $ as an invariant rather than $\phi $.
This is illusory. The problem is that the map $\phi _K $ is only defined if we understand
$End _{\ZZ /2}  (K^*(X))$ to mean {\emph{endomorphisms of $\ZZ /2$-graded abelian groups}} since
 if $A$ is non-commutative then {\emph{there is no ring structure}} on $K_*(A)$.

 We can say something about the map $\phi $. If the map $h: A \to A$ arises from
  a map $f: X \to X$ then of course $\phi (h) = f^* \otimes 1 $.  If $h $ arises from conjugation
 by a unitary then $\phi  (h)$ is just the identity, since $U_n$ is path-connected.    We hope to compute $\phi (h)$
 in some intertwined cases and expect it to be a helpful invariant. This is work in progress.

%%%%%%%%%%%%%%%%%%%%%%%%%%%%%%%%%%%%%%%%%%%%%%%%%%%%%%
%%%%%%%%%%%%%%%%%%%%%%%%%%%%%%%%%%%%%%%%%%%%%%%%%%%%%
%%%%%%%%%%%%%%%%%%%%%%%%%%%%%%%%%%%%%%%%%%%%%%%%%%%%%%%%%%

%%%%%%%%%%%%%%%%%%%%%%%%%%%%%%%%%%%%%%%%%%%%%%%%%%%%%%%%
%%%%%%%%%%%%%%%%%%%%%%%%%%%%%%%%%%%%%%%%%%%%%%%%%%%%%%%%%

%%%%%%%%%%%%%%%%%%%%%%%%%%%%%%%%%%%%%%%%%%%%%%%%%%%%%
%%%%%%%%%%%%%%%%%%%%%%%%%%%%%%%%%%%%%%%%%%%%%%%%%%%%%%%%
%%%%%%%%%%%%%%%%%%%%%%%%%%%%%%%%%%%%%%%%%%%%%%%%%%%%%%%55
%%%%%%%%%%%%%%%%%%%%%%%%%%%%%%%%%%%%%%%%%%%%%%%%%%%%%%5
%%%%%%%%%%%%%%%%%%%%%%%%%%%%%%%%%%%%%%%%%%%%%%%%%%%%%%%%
%%%%%%%%%%%%%%%%%%%%%%%%%%%%%%%%%%%%%%%%%%%%%%%%%%%%%%%%%%
%%%%%%%%%%%%%%%%%%%%%%%%%%%%%%%%%%%%%%%%%%%%%%%%%%%%%%%

\end{document}